# THE FIXED POINT PROPERTY
# FOR POSETS OF SMALL WIDTH

Jonathan David Farley

Abstract. The fixed point property for finite posets of width 3 and 4 is studied in terms of forbidden retracts. The ranked forbidden retracts for width 3 and 4 are determined explicitly. The ranked forbidden retracts for the width 3 case that are linearly indecomposable are examined to see which are minimal automorphic. Part of a problem of Niederle from 1989 is thus solved.

## 1. Introduction.

How does one recognize a poset with the fixed point property? The problem has a long history ([12]). Recently, it was shown to be computationally intractable ([4], Theorem 1.1). Nonetheless, for certain classes of ordered sets, the posets with the fixed point property *have* been characterized, generally in terms of forbidden retracts:

A finite, connected rank 1 poset has the fixed point property if and only if no crown is a retract. (There is also a characterization for the infinite rank 1 posets in terms of forbidden suborders. See [11], Theorem; [10], T heorem 4.) A finite, $N$-free poset of dimension 2 has the fixed point property if and only if it has no 4-cycle tower retract ([6], Theorem 4). A finite width 2 poset has the fixed point property if and only if no 4-crown tower of t he poset is a retract. (More generally, the chain-complete posets have been characterized. See [7], Theorem 1.) A finite width 3 poset has the fixed point property if and only if no tower of sections is a retract (if and only if no tow er of nice or very nice sections is a retract). These sections, even the "very nice" ones, have a rather complicated description ([9], Theorem). Width 3 posets containing an element

1991 *Mathematics Subject Classification*. 06A06.

The author would like to thank Prof. Bernd Schröder for pointing out complexities in the non-ranked case, which the author had overlooked.

Part of this research was conducted at the Mathematical Sciences Research Institute.

Typeset by $\mathcal{AMS}$-TEX

comparable to every other element are studied in [ 8], Theorems 1 and 2. A finite poset has the fixed point property if and only if no generalized crown is a retract ([1], Theorem 1).

We continue the quest: We explicitly describe the ranked forbidden retracts for finite width 4 posets with the fixed point property (Theorem 3.5). As a corollary, we obtain a very simple description of the ranked forbidden retracts of width 3 posets (cf. [9], Theorem): See Corollary 3.6, Proposition 4.1, and Proposition 4.2. (Results of B. S. W. Schröder suggest that the non-ranked case may be intractable [13].)

The posets in our list of forbidden retracts are ordinal sums of ranked posets. The elements in consecutive ranks therefore form bipartite posets with at most four minimal elements and four maximal elements. We explicitly describe these bipartite posets .

Finally, we solve part of a two-part problem proposed by Niederle in 1989 ([9]): to see if the linearly indecomposable forbidden retracts that he found for the width 3 case are minimal, in the sense that they do not admit proper retracts of a similar form (Theorem 5.12).

## 2. Notation and Definitions.

*All posets are finite.*[1] For basic terminology, consult [3].

Let $P$ be a poset. Its cardinality will be denoted $\#P$.

A function $f : P \to P$ is *order-preserving* if, for all $p, q \in P$, $p \leqslant q$ implies $f(p) \leqslant f(q)$. The function is an *automorphism* if it is bijective and both it and its inverse are order-preserving. A subset $Q \subseteq P$ is a *retract* if there exists an order-preserving map $f : P \to Q$ such that $f(q) = q$ for all $q \in Q$; the map $f$ is a *retraction.*

An element $p \in P$ is a *fixed point* of a function $f : P \to P$ if $p = f(p)$. An order-preserving map with no fixed points is *fixed-point-free.* The poset $P$ has the *fixed point property* if every order-preserving map on $P$ has a fixe d point. It is well known that a poset has the fixed point property if and only if every retract does ([5], §1).

A poset is *automorphic* if it has a fixed-point-free automorphism; an automorphic poset is *minimal* if no proper retract is automorphic ([2]).

Let $P$ and $Q$ be posets: $P + Q$ is the poset with underlying set $P \cup Q$ such that no point of $P$ is comparable with a point of $Q$. The *ordinal sum* of $P$ and $Q$ is the poset with underlying set $P \cup Q$ such that $p < q$ for all $p \in P$, $q \in Q$. (Similarly, one may define ordinal sums of zero, one, or more than two posets.)

An *antichain* is a poset in which no two distinct elements are comparable. The *width* of a poset is the cardinality of its largest antichain.

---

[1] This statement is known as Trotter's Axiom.

A *chain* is a totally ordered set $C$; its *rank* is $^\#C - 1$. The *height* of a poset $P$ is the rank of its longest chain. A poset $P$ is *ranked* of rank $\mathrm{r}(P) \in \mathbb{N}_0$ if every maximal chain has rank $\mathrm{r}(P)$. The *rank* of an element $p \in P$ is the height $\mathrm{r}(p)$ of the poset $\{\, q \in P \mid q \leqslant p \,\}$. If $i \leqslant j$ ($i, j \in \mathbb{N}_0$), then $P(i, j)$ is the poset $\{\, p \in P \mid i \leqslant \mathrm{r}(p) \leqslant j \,\}$; let $P(i) := P(i, i)$.

Let $P$ be a poset and $p, q \in P$. We say $p$ is a *lower cover* of $q$ (and $q$ is an *upper cover* of $p$) if $p < q$ and $p \leqslant r < q$ implies $p = r$ for all $r \in P$. An element $p \in P$ is *irreducible* if it has a unique upper cover or a unique lower cover.

For $n \geqslant 2$, the *$2n$-crown* $C_{2n}$ is the poset with underlying set

$$\{\, x_0, \ldots, x_{n-1}; y_0, \ldots, y_{n-1} \,\}$$

in which

$$x_0 < y_0 > x_1 < y_1 > \cdots > x_{n-1} < y_{n-1} > x_0$$

are the only strict comparabilities. Figure 2.1 shows $C_4$, $C_6$, and $C_8$.

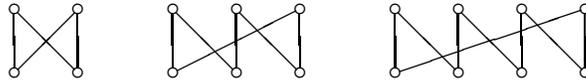

**Figure 2.1.** The crowns $C_4$, $C_6$, and $C_8$.

A *4-tower* is an ordinal sum of two-element antichains (Figure 2.2).

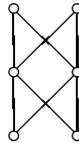

**Figure 2.2.** An example of a 4-tower.

[In [7], a *4-crown tower* in a poset $P$ is defined to be a 4-tower $T$ such that each $t \in T$ is minimal (in $P$) in the set of upper bounds of $\{\, t' \in T \mid t' < t \,\}$. In [6], a *4-cycle tower* in a poset $P$ is defined to be a 4-tower $T$ such that, for $0 \leqslant i < \mathrm{r}(T)$, no element of $P$ is an upper bound of $T(i)$ and a lower bound of $T(i+1)$.]

A *6-stack* is a poset of rank $n \geqslant 1$ such that, for $0 \leqslant i < n$, $P(i, i+1)$ is a 6-crown (Figure 2.3).

A *6-tower* is an ordinal sum of two-element antichains and 6-stacks (Figure 2.4).

An *8-stack* is a poset $P$ of rank $n \geqslant 1$ such that, for $0 \leqslant i < n$, $P(i, i+1)$ is isomorphic to one of the following posets with four minimal and four maximal elements:

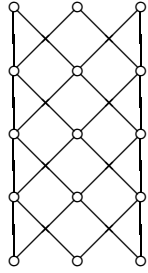

**Figure 2.3.** An example of a 6-stack.

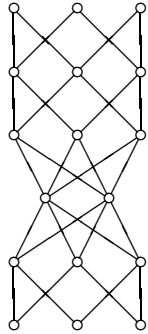

**Figure 2.4.** An example of a 6-tower.

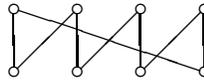

**Figure 2.5a.** The crown $C_8$.

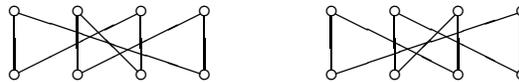

**Figure 2.5b.** The crown $C_8$.

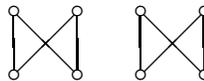

**Figure 2.6a.** The poset $C_4 + C_4$.

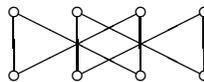

**Figure 2.6b.** The poset $C_4 + C_4$.

(1) $C_8$ (Figures 2.5a and 2.5b);
(2) $C_4 + C_4$ (Figures 2.6a–2.6c);

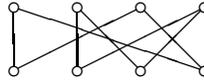

**Figure 2.6c.** The poset $C_4 + C_4$.

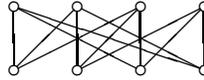

**Figure 2.7a.** The poset $\bar{K}_{4,4}$.

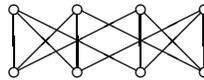

**Figure 2.7b.** The poset $\bar{K}_{4,4}$.

(3) the almost-complete bipartite poset $\bar{K}_{4,4}$ (Figures 2.7a and 2.7b);
(4) the poset $Z_1$ or its dual $Z_1^\partial$ (Figures 2.8a–2.8b and 2.9);

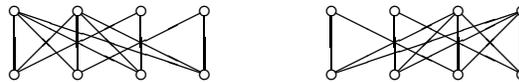

**Figure 2.8a.** The poset $Z_1$.

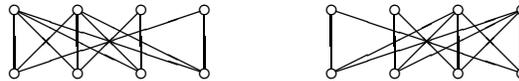

**Figure 2.8b.** The poset $Z_1$.

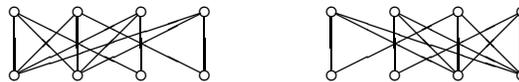

**Figure 2.9.** The poset $Z_1^\partial$.

(5) the poset $Z_2$ (Figures 2.10a and 2.10b);

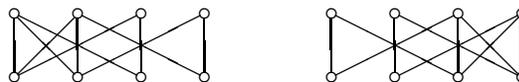

**Figure 2.10a.** The self-dual poset $Z_2$.

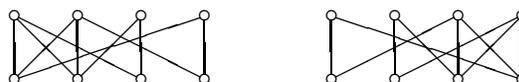

**Figure 2.10b.** The self-dual poset $Z_2$.

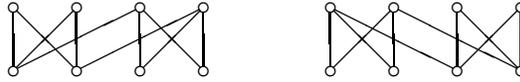

**Figure 2.11.** The self-dual poset $Z_3$.

(6) the poset $Z_3$ (Figure 2.11);
(7) the poset $Z_4$ (Figure 2.12);

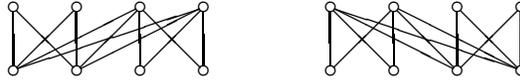

**Figure 2.12.** The self-dual poset $Z_4$.

(8) the poset $Z_5$ (Figure 2.13);

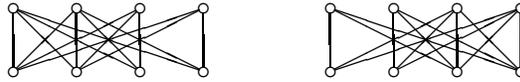

**Figure 2.13.** The self-dual poset $Z_5$.

The 8-stack is *admissible* if there exist $F\colon \{0,\ldots,n\} \to \{(4),(2)(2)\}$ and a *drawing* of $P$ such that the following holds. For $0 \leqslant i < n$, the portion of the drawing corresponding to $P(i, i+1)$ looks like one of the following figure s, according to the value of $F(i)$ and $F(i+1)$ (Table 2.1):

| $F(i)$ vs. **F(i+1)** | **(4)** | **(2)(2)** |
|---|---|---|
| (4) | 2.5a, 2.6b, 2.7a | 2.6b, 2.8a |
|  | $C_8$, $C_4 + C_4$, $\bar{K}_{4,4}$ | $C_4 + C_4$, $Z_1$ |
|  |  |  |
| (2)(2) | 2.6b, 2.6c, 2.9 | 2.5b, 2.6a, 2.6b, 2.6c, 2.7b, 2.8a, 2.8b, 2.9, 2.10a, 2.10b, 2.11, 2.12, 2.13 |
|  | $C_4 + C_4$, $C_4 + C_4$, $Z_1^\partial$ | $C_8$, $C_4 + C_4$, $C_4 + C_4$, $C_4 + C_4$, $\bar{K}_{4,4}$, $Z_1$, $Z_1$, $Z_1^\partial$, $Z_2$, $Z_2$, $Z_3$, $Z_4$, $Z_5$ |

**Table 2.1.** The figures that resemble $P(i, i+1)$ in the drawing of $P$.

An *8-tower* is an ordinal sum of 6-towers and admissible 8-stacks.

Let $m, n, k \in \mathbb{N}_0$ be such that $m, n \geqslant 2$ and $0 \leqslant k \leqslant n$. Let $j_1, \ldots, j_k \in \mathbb{N}_0$ be such that $0 \leqslant j_1 < j_2 < \cdots < j_k \leqslant n-1$. Consider the poset with underlying set $\{x_0, \ldots, x_{m-1}, y_0, \ldots, y_{n-1}\}$ in which

$$x_i < y_{i+j_1}, \ldots, y_{i+j_k} \qquad (i = 0, \ldots, n-1)$$

are the only strict comparabilities. (Calculate the subscripts modulo $n$.) The poset is *circulant* if the function sending $x_i$ to $x_{i+1}$ ($i = 0, \ldots, m-1$) and $y_j$ to $y_{j+1}$ ($j = 0, \ldots, n-1$) is order-preserving. In this case, the pair

$$\{(x_0, \ldots, x_{m-1}), (y_0, \ldots, y_{n-1})\}$$

is the *bipartition* of the poset.

Suppose that, for some $k \in \mathbb{N}_0$, the poset $P$ can be partitioned into cyclic sequences $A_0, \ldots, A_{k-1}$, each an antichain with at least two elements, such that, for $0 \leqslant i < j \leqslant k-1$, the subposet $A_i \cup A_j$ is circulant with bipartition $\{A_i, A_j\}$. Then $P$ is called a *generalized crown* with *k-partition* $\{A_0, \ldots, A_{k-1}\}$. (The definitions of this and the previous paragraph come from [1].)

A generalized crown with a $k$-partition ($k \in \mathbb{N}_0$) is *special* if it admits an automorphism whose orbits are the members of the $k$-partition.

A *section* is either a two-element antichain or else a poset with underlying set $\{[i, k] \mid 0 \leqslant i \leqslant 2;\ 0 \leqslant k \leqslant n\}$ where $n \geqslant 1$ and
  (1) if $0 \leqslant i \leqslant 2$ and $0 \leqslant k < l \leqslant n$, then $[i, k] < [i, l]$;
  (2) if $0 \leqslant k \leqslant n$, then $\{[0, k], [1, k], [2, k]\}$ is an antichain;
  (3) if $0 \leqslant i, j \leqslant 2$ and $0 \leqslant k, l \leqslant n$, then $[i, k] < [j, l]$ implies $[i+1, k] < [j+1, l]$ (with addition modulo 3);
  (4) if $0 \leqslant k < n$, then for some $i, j \in \{0, 1, 2\}$, $[i, k] \not< [j, k+1]$.

A section $P$ is *nice* if, for all $p, q \in P$ such that $p < q$, there exist $r, s \in P$ such that $p < r$ but $q \not< r$, and $s < q$ but $s \not< p$. A *tower of sections* is an ordinal sum of sections. A section is *very nice* if no pr oper retract is isomorphic to a tower of sections. (These definitions come from [9].)

## 3. The Fixed Point Property for Ranked Posets of Width 4.

In this section, we describe the ranked minimal automorphic posets of width at most 4 (Theorem 3.5). As a corollary, we obtain a characterization of the width 3 posets (Corollary 3.6) (cf. [9], Theorem).

The proof of [1], Theorem 1 actually establishes the following lemma:

**Lemma 3.1.** *A poset without the fixed point property has a special generalized crown as a retract.*   □

The next lemma is trivial.

**Lemma 3.2.** *Let $\{A_0, \ldots, A_{k-1}\}$ be the $k$-partition of a special generalized crown ($k \in \mathbb{N}_0$). If $r \in \mathbb{N}_0$ and $0 \leqslant i_1 < \cdots < i_r \leqslant k-1$, then $\{A_{i_1}, \ldots, A_{i_r}\}$ is an $r$-partition of a special generalized crown.*   □

Schröder has proven the following fact [13]:

**Fact 3.3.** *There exists a width 3 minimal automorphic special generalized crown that is not ranked.*

Finally, we have:

**Lemma 3.4.** *Every special generalized crown has a special generalized crown with no irreducibles as a retract.*

*Proof.* If an element is irreducible, so is every element in its orbit. By dismantling by irreducibles (see [14], Theorem 2.6) and Lemma 3.2, we get the result. □

**Theorem 3.5.** *A (finite) ranked poset of width at most 4 that is minimal automorphic is an 8-tower; and every 8-tower is automorphic.*

*Proof.* Every admissible 8-stack $P$ is automorphic. Define a fixed-point-free automorphism $f : P \to P$ as follows: Consider an admissible diagram. Let $0 \leqslant i \leqslant r(P)$. If $F(i) = (4)$, let $f$ cycle the rank $i$ elements from left to right; if $F(i) = (2)(2)$, let $f$ switch the two leftmost elements as well as the two rightmost elements.

Now assume we have a ranked minimal automorphic poset $R$ of width at most 4. By Lemmas 3.1 and 3.4, it is a special generalized crown with no irreducibles. Let $f: R \to R$ be the fixed-point-free automorphism.

Since $R$ is ranked, we can draw $R$ so that the elements from each rank are side-by-side. There are at most 4 elements in each rank. Without loss of generality, $f$ either cyclically permutes the elements from left to right in the diagram, or else it switches the two leftmost elements as well as the two rightmost elements. Accordingly, set $F(i)$ equal to (2), (3), (4), or (2)(2), where $0 \leqslant i \leqslant r(R)$.

The possible diagrams of $R(i, i+1)$ may be determined by combining the above four types, bearing in mind that no element is irreducible, and using the action of $f$; we write $\alpha/\beta$ if $F(i) = \alpha$ and $F(i+1) = \beta$ $[0 \leqslant i < r(R)]$. We allow ourselves to manipulate the maximal elements but not the minimal ones.

For (2)/(2), (2)/(3), (2)/(4), or (2)/(2)(2), we get an ordinal sum (Figures 3.1–3.3).

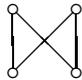

**Figure 3.1.** Case (2)/(2).

For (3)/(3), there are two possibilities: $C_6$ or an ordinal sum (Figure 3.4).
For (3)/(4) and (3)/(2)(2), one gets an ordinal sum (Figure 3.5).
For (4)/(4), one gets $C_8$, $C_4 + C_4$, $\bar{K}_{4,4}$, or an ordinal sum (Figure 3.6).

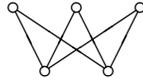

**Figure 3.2.** Case (2)/(3).

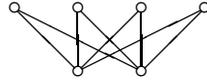

**Figure 3.3.** Cases (2)/(4) and (2)/(2)(2).

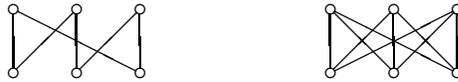

**Figure 3.4.** Case (3)/(3).

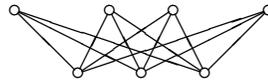

**Figure 3.5.** Cases (3)/(4) and (3)/(2)(2).

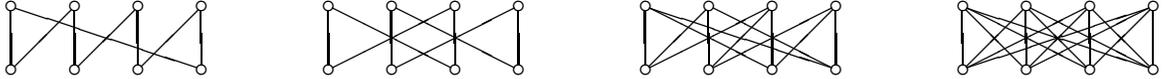

**Figure 3.6.** Case (4)/(4).

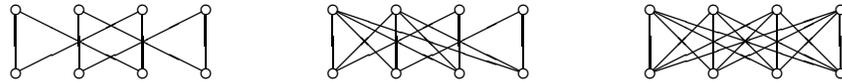

**Figure 3.7.** Case (4)/(2)(2).

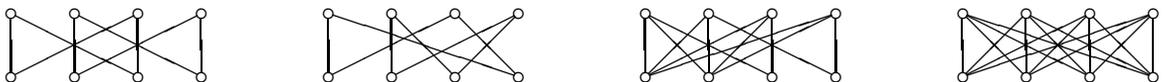

**Figure 3.8.** Case (2)(2)/(4).

For (4)/(2)(2), one gets $C_4 + C_4$, $Z_1$, or an ordinal sum (Figure 3.7).

For (2)(2)/(4), one gets $C_4 + C_4$, $Z_1^\partial$, or an ordinal sum (Figure 3.8).

Finally, consider (2)(2)/(2)(2). We use the notation $\langle a, b \rangle$ to indicate the number of upper covers the four minimal elements have ($2 \leqslant a, b \leqslant 4$).

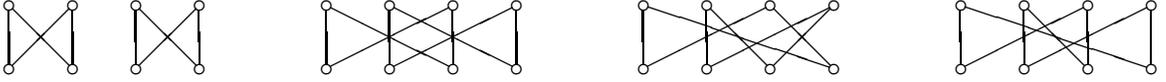

**Figure 3.9.** Case $(2)(2)/(2)(2);\langle 2,2\rangle$.

For $\langle 2,2\rangle$, one gets $C_4 + C_4$ or $C_8$ (Figure 3.9).
For $\langle 2,3\rangle$, we get $Z_2$ or $Z_3$ (Figure 3.10).

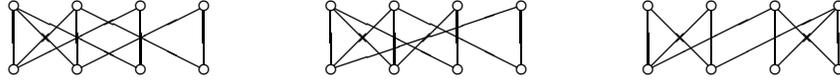

**Figure 3.10.** Case $(2)(2)/(2)(2);\langle 2,3\rangle$.

For $\langle 2,4\rangle$, we get $Z_1^\partial$ or $Z_4$ (Figure 3.11).

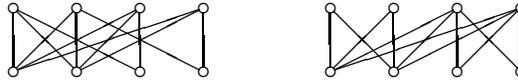

**Figure 3.11.** Case $(2)(2)/(2)(2);\langle 2,4\rangle$.

For $\langle 3,3\rangle$, we get $\bar{K}_{4,4}$ or $Z_1$ (Figure 3.12).

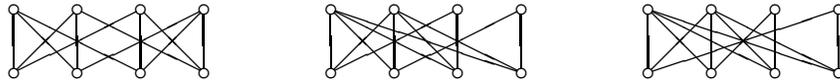

**Figure 3.12.** Case $(2)(2)/(2)(2);\langle 3,3\rangle$.

For $\langle 3,4\rangle$, we get $Z_5$ (Figure 3.13).

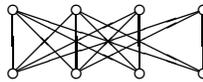

**Figure 3.13.** Case $(2)(2)/(2)(2);\langle 3,4\rangle$.

For $\langle 4,4\rangle$, we get an ordinal sum (Figure 3.14).

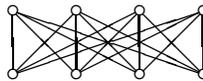

**Figure 3.14.** Case $(2)(2)/(2)(2);\langle 4,4\rangle$.

Three- or four-element antichains that are ordinal summands can be retracted onto two-element antichains.

The theorem is proved. □

We get the following corollary for width 3 posets:

**Corollary 3.6.** *A (finite) ranked poset of width at most 3 that is minimal automorphic is a 6-tower; and every 6-tower is automorphic.* □

## 4. The Structure of the Ranked Forbidden Retracts in the Width 3 Case.

In §3, width 3 posets with the fixed point property were described in terms of forbidden retracts, the *6-towers*. The 6-towers are easily described; in this section, we show that they can be even more easily visualized (Proposition 4.1 and Figure 4.2). In general, the width 3 forbidden retracts introduced in [9], the towers of (nice or very nice) sections, have a fairly complicated definition, and it is difficult to see what these posets look like from that definition. Nonetheless, we show that the ranked towers of nice sections are precisely our 6-towers (Proposition 4.2).

The first proposition improves our description of width 3 posets with the fixed point property (Corollary 3.6). A priori, there are many different ways to stack 6-crowns (Figure 4.1). In reality, there is only one.

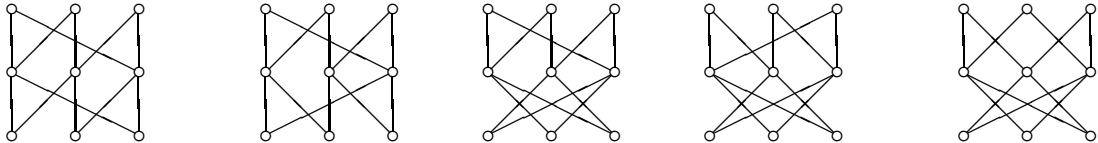

**Figure 4.1.** Various ways of stacking 6-crowns.

**Proposition 4.1.** *For each $n \geqslant 1$, there is (up to isomorphism) only one 6-stack of rank $n$ (Figure 4.2).*

*Proof.* Given a 6-stack, draw it so that the elements of each rank form a row. Let $x_0$, $y_0$, and $z_0$ be the minimal elements from left to right. Shift the rank 1 elements so that the leftmost element covers $x_0$ and $y_0$, the middle element $x_0$ and $z_0$, and the rightmost element $y_0$ and $z_0$. Proceed to the next rank. □

We thus see that every 6-stack is a *superposition of crowns* as defined in [2].

Now we give a nice description of the rather cumbersome-to-define nice sections of [9] (at least when they are ranked):

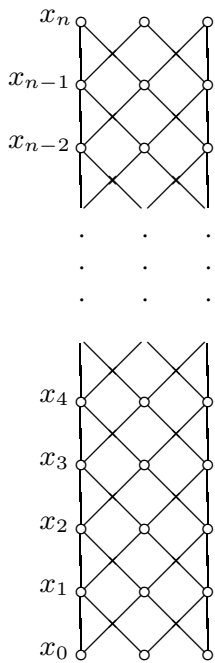

**Figure 4.2.** The unique 6-stack of rank $n \geqslant 1$.

**Proposition 4.2.** *The ranked nice sections are exactly the two-element antichains and the 6-stacks. Hence, every ranked tower of nice sections is a 6-tower, and conversely.*

*Proof.* It is easy to see (especially in light of Proposition 4.1) that every 6-stack is a nice section.

Let $P$ be a ranked nice section that is not an antichain. Conditions (1) and (2) of the definition show that every rank has three elements. Conditions (3) and (4) say that no element has three upper covers or three lower covers. Since $P$ is nice, no element is irreducible.

Hence, every non-maximal element has exactly two upper covers, and every non-minimal element has exactly two lower covers. That is, $P(i, i+1)$ is a 6-crown for $0 \leqslant i < \mathrm{r}(P)$. □

The examples of B. S. W. Schröder show that, in the non-ranked setting, even nice sections can be nasty [13].

## 5. Nice Sections That Are Very Nice: A Problem of Niederle.

In [9], the following question appears:

**Problem** (Niederle, 1989). *Are there nice sections that are not very nice? Characterize very nice sections.*

We determine which ranked nice sections are very nice in Theorem 5.12; together with Propositions 4.1 and 4.2, this solves part of Niederle's problem.

**Lemma 5.1.** *Let $P$ be a 6-stack and $Q \subseteq P$ a tower of sections. If $^\#Q(j) = 3$ (where $j \in \mathbb{N}_0$), then $Q(j) = P(i)$ for some $i \in \mathbb{N}_0$.*

*Proof.* The lemma follows by noting that, if $x$, $y$, $z \in P$ form an antichain, they must have the same rank. (See Figure 4.2.) □

**Lemma 5.2.** *Let $P$ be a 6-stack and $Q \subseteq P$ a tower of sections. If $P(i) \subseteq Q$ for some $i \in \mathbb{N}_0$ but $P(i+1) \not\subseteq Q$, then $P(i+1) \subseteq P \setminus Q$.*

*Proof.* Assume for a contradiction that $P(i+1) \cap Q(j) \neq \emptyset$ for some $j \in \mathbb{N}_0$. By Lemma 5.1, $Q(j)$ must be a two-element ordinal summand of $Q$, so that an element of $P(i+1)$ has three lower covers, which is absurd. □

**Lemma 5.3.** *Let $P$ be a 6-stack, $i \in \mathbb{N}_0$, and let $p \in P(i)$ and $q \in P(i+2)$. Then $p < q$.*

*Proof.* See Figure 4.2. □

**Lemma 5.4.** *Let $P$ be a 6-stack and $Q$ a retract that is a tower of sections. If $\emptyset \neq P(i) \subseteq Q$ for some $i \in \mathbb{N}_0$, then $P = Q$.*

*Proof.* Suppose not, for a contradiction. Then, without loss of generality, $P(i) \subseteq Q$ but $P(i+1) \not\subseteq Q$. By Lemma 5.2, $P(i+1) \subseteq P \setminus Q$. Because $Q$ is a retract, there must be some $j \geq i+2$ such that $P(j) \cap Q \neq emptyset$; choose $j$ minimal and let $k \in \mathbb{N}_0$ be such that $P(j) \cap Q(k) \neq \emptyset$. Thus some element of $P(i+1)$ is a lower bound of two elements in $Q(k)$, so it has nowhere to go under a retraction. □

**Corollary 5.5.** *If a 6-stack has a proper retract that is a tower of sections, then that retract must be a 4-tower.*

*Proof.* The corollary follows from Lemmas 5.1 and 5.4. □

In the sequel, we will consider a 6-stack $P$ of rank $n \geq 1$, drawn in the fashion prescribed by the proof of Proposition 4.1. Label the elements of $P(i)$ from left to right $x_i$, $y_i$, and $z_i$ ($0 \leq i \leq n$). See Figure 5.1.

**Lemma 5.6.** *Let $P$ be a 6-stack and $Q$ a 4-tower retract. Then $P(0) \cap Q(0) \neq \emptyset$.*

*Proof.* If $P(0) \cap Q(0) = \emptyset$, then $P(0)$ contains a lower bound of $Q(0)$, which thus has nowhere to go under a retraction. □

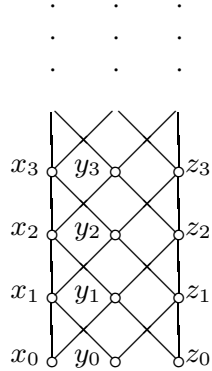

**Figure 5.1.** The 6-stack $P$.

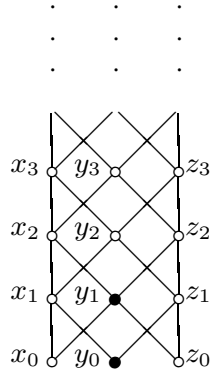

**Figure 5.2.** The set $Q(0)$ in $P$.

**Lemma 5.7.** *Let $P$ be a 6-stack and $Q$ a 4-tower retract such that $^\#P(0) \cap Q(0) = 1$.*

*Then $n \geqslant 3$ and $P(3,n)$ has a 4-tower retract.*

*Proof.* Without loss of generality, $Q(0) = \{y_0, y_1\}$ (Figure 5.2).

Let $f \colon P \to Q$ be the retraction. Then

$$f(x_0) = y_1 = f(z_0),$$

so $f(x_1), f(z_1) \in P(2,n)$.

$$\text{Case 1. } {}^\#P(2) \cap Q(1) = 2$$

Then $Q(1) = \{x_2, z_2\}$ (Figure 5.3).

Hence $f(x_1) = x_2$, $f(z_1) = z_2$, and $f(y_2), f(y_3) \in P(3,n)$. Therefore, $f$ preserves $P(3,n)$ and retracts it onto a 4-tower.

$$\text{Case 2. } {}^\#P(2) \cap Q(1) = 1$$

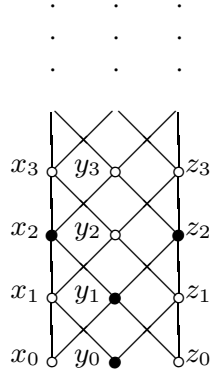

**Figure 5.3.** The set $Q(0,1)$ in $P$.

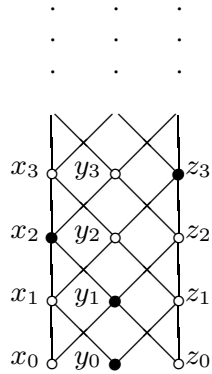

**Figure 5.4.** The set $Q(0,1)$ in $P$.

Without loss of generality, $Q(1) = \{x_2, z_3\}$ (Figure 5.4).
Hence $x_2 = f(x_1) \leqslant f(z_3) = z_3$, a contradiction.

$$\text{Case 3. } P(2) \cap Q(1) = \emptyset$$

Then $f(x_1) \in P(3,n)$, so $f$ preserves $P(3,n)$ and retracts it onto a 4-tower. $\square$

**Lemma 5.8.** *Let $P$ be a 6-stack and $Q$ a 4-tower retract such that $^{\#}P(0) \cap Q(0) = 2$.*
*Then $n \geqslant 3$ and $P(3,n)$ has a 4-tower retract.*

*Proof.* Without loss of generality, $Q(0) = \{x_0, z_0\}$ (Figure 5.5). Let $f : P \to Q$ be the retraction.

$$\text{Case 1. } P(1) \cap Q(1) \neq \emptyset$$

Then $Q(1) = \{y_1, y_2\}$ (Figure 5.6).

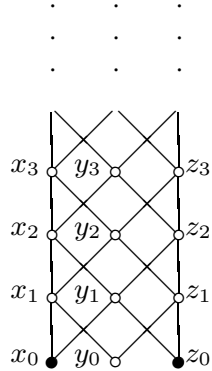

**Figure 5.5.** The set $Q(0)$ in $P$.

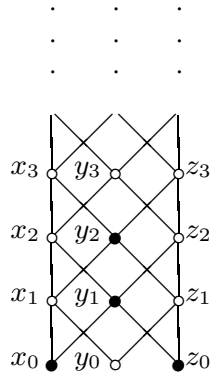

**Figure 5.6.** The set $Q(0,1)$ in $P$.

We are done if $f(x_2), f(z_2) \in P(3,n)$, so we may assume that $f(x_2) = y_1$. Then $f(x_1) = x_0 = f(y_0)$ and $f(z_1) = y_2$. Hence $f(z_2) \in P(3,n)$. Also, $f(x_3) \in P(3,n)$, so $f$ preserves $P(3,n)$ and retracts it onto a 4-tower.

*Case 2.* $P(1) \cap Q(1) = \emptyset$ and $^\#P(2) \cap Q(1) = 2$

Because of $y_1$, we must have (without loss of generality) $Q(1) = \{x_2, y_2\}$ (Figure 5.7).

Hence $f(x_1) = x_0 = f(y_0)$, $f(y_1) = x_2$, and $f(z_1) = y_2$. Therefore, $f(z_2) \geqslant x_2, y_2$, so $n \geqslant 3$, and $f(x_3) \geqslant x_2, y_2$, so $f$ preserves $P(3,n)$ and retracts it onto a 4-tower.

*Case 3.* $P(1) \cap Q(1) = \emptyset$ and $^\#P(2) \cap Q(1) = 1$

Because of $y_1$, we must have $Q(1) = \{y_2, y_3\}$ (Figure 5.8).
Thus $f(x_1) = x_0$, $f(z_1) = z_0$, and $y_0$ has nowhere to go under $f$.

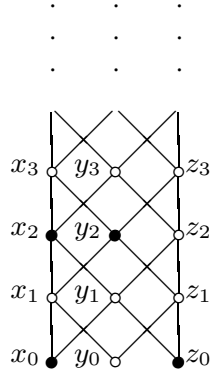

**Figure 5.7.** The set $Q(0,1)$ in $P$.

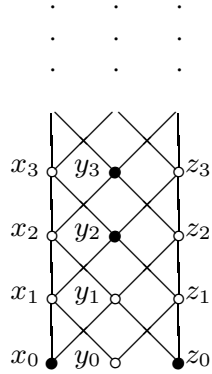

**Figure 5.8.** The set $Q(0,1)$ in $P$.

*Case 4.* $P(1) \cap Q(1) = \emptyset$ and $P(2) \cap Q(1) = \emptyset$

Now $y_1$ has nowhere to go under $f$. $\square$

**Lemma 5.9.** *Let $P$ be the 6-stack of rank 3 in Figure 5.9. Let $Q$ be the 4-tower $\{x_0, z_0, y_1, y_2, x_3, z_3\}$.*

*Then there is a unique retraction $f\colon P \to Q$ such that $f(x_1) = x_0$; under this retraction, $f(y_0) = x_0$ and $f(y_3) = z_3$.*

*Proof.* If $f\colon P \to Q$ is a retraction such that $f(x_1) = x_0$, then $f(y_0) = x_0$, $f(z_1) = y_2$, $f(z_2) = z_3$, $f(y_3) = z_3$, and $f(x_2) = y_1$.

Conversely, the map $f : P \to Q$ defined as above (and the identity on $Q$) is a retraction (Figure 5.10). $\square$

The following is an easy corollary.

**Corollary 5.10.** *A 6-stack whose rank is a multiple of 3 has a 4-tower as a retract.* $\square$

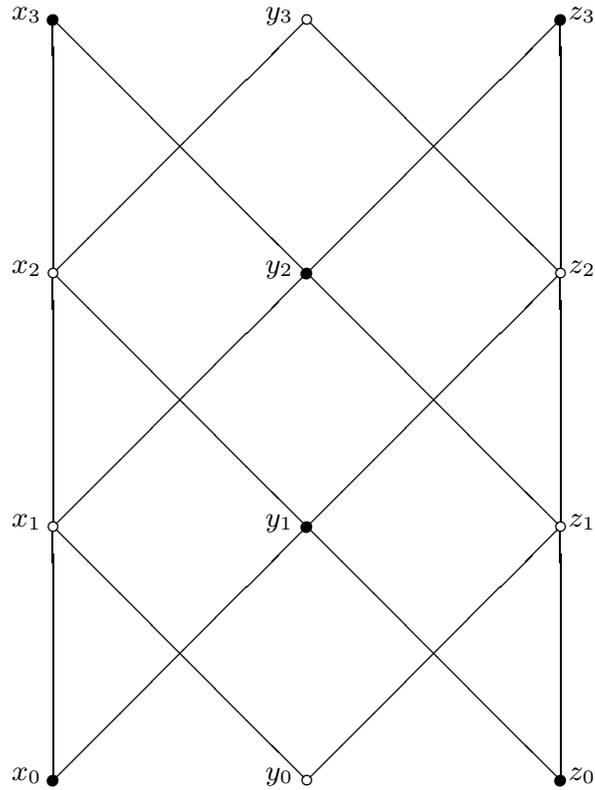

**Figure 5.9.** The 6-stack $P$ and 4-tower $Q$.

**Proposition 5.11.** *Let $P$ be a 6-stack.*

*Then $P$ has a proper retract that is a tower of sections if and only if $\mathrm{r}(P)$ is a multiple of 3.*

*Therefore, $P$ is minimal automorphic if and only if $\mathrm{r}(P)$ is not a multiple of 3.*

*Proof.* Assume that $P$ has a proper retract $Q$ that is a tower of sections. By Corollary 5.5, $Q$ is a 4-tower. By Lemmas 5.6–5.8, $n \geqslant 3$ and $P(3,n)$ has a 4-tower retract. Continuing, we get a contradiction unless $n$ is a multiple of 3. The converse follows from Corollary 5.10.

The final statement comes from [9], Theorem. □

**Theorem 5.12.** *A ranked nice section is very nice if and only if its rank is not a positive multiple of 3.*

*Proof.* The theorem follows from Propositions 4.2 and 5.11. □

Part of the 1989 problem of Niederle ([9]) is thus solved.

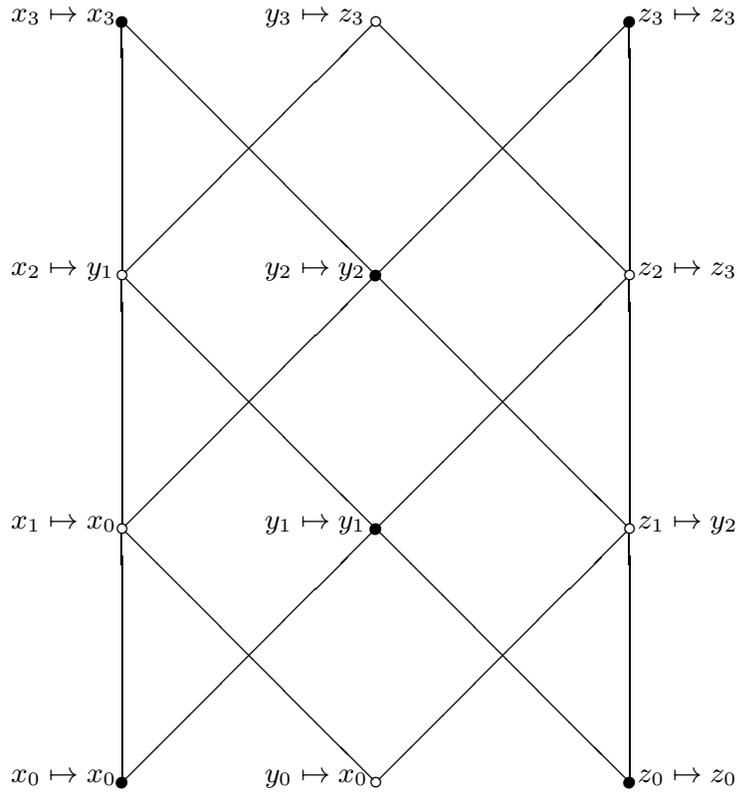

**Figure 5.10.** The retraction $f\colon P \to Q$.

Department of Mathematics, Vanderbilt University, Nashville, Tennessee 37240, United States of America

Mathematical Sciences Research Institute, 1000 Centennial Drive, Berkeley, California 94720, United States of America

farley@math.vanderbilt.edu